\documentclass{article}
\usepackage{amssymb,amsmath,amsfonts,amsthm}
\usepackage{pgf,tikz,xcolor}
\usepackage{pgfplots}
\usetikzlibrary{arrows}
\usepackage{graphicx}
\usepackage{subfig}
\usepackage[colorlinks=true,citecolor=blue,urlcolor=blue,linkcolor=blue]{hyperref}
\usepackage{lineno}
 \newcommand*\patchAmsMathEnvironmentForLineno[1]{%
  \expandafter\let\csname old#1\expandafter\endcsname\csname #1\endcsname
  \expandafter\let\csname oldend#1\expandafter\endcsname\csname end#1\endcsname
  \renewenvironment{#1}%
     {\linenomath\csname old#1\endcsname}%
     {\csname oldend#1\endcsname\endlinenomath}}% 
\newcommand*\patchBothAmsMathEnvironmentsForLineno[1]{%
  \patchAmsMathEnvironmentForLineno{#1}%
  \patchAmsMathEnvironmentForLineno{#1*}}%
\AtBeginDocument{%
\patchBothAmsMathEnvironmentsForLineno{equation}%
\patchBothAmsMathEnvironmentsForLineno{align}%
\patchBothAmsMathEnvironmentsForLineno{flalign}%
\patchBothAmsMathEnvironmentsForLineno{alignat}%
\patchBothAmsMathEnvironmentsForLineno{gather}%
\patchBothAmsMathEnvironmentsForLineno{multline}%
}

\newcommand{\dt}{\Delta t}

\usepackage[ruled,vlined]{algorithm2e}
\RestyleAlgo{boxed}
\newcommand{\Dt}{\Delta t}

\title{Revisionist Integral Deferred Correction with Adaptive Step-Size
  Control}

\author{A.~J.~Christlieb \and C.~B.~Macdonald \and
  B.~W.~Ong\footnote{corresponding author, ongbw@msu.edu} \and R.~J.~Spiteri}

\begin{document}
%\linenumbers

\maketitle
\begin{abstract}
  Adaptive step-size control is a critical feature for the robust and
  efficient numerical solution of initial-value problems in ordinary
  differential equations. In this paper, we show that adaptive
  step-size control can be incorporated within a family of parallel
  time integrators known as Revisionist Integral Deferred Correction
  (RIDC) methods.  The RIDC framework allows for various strategies to
  implement step-size control, and we report results from exploring a
  few of them.  \bigskip


  %\bigskip
  %\noindent {\em AMS Classifications}:
 \end{abstract}

\section{Introduction}

The purpose of this paper is to show that local error estimation and
adaptive step-size control can be incorporated in an effective manner
within a family of parallel time integrators based on Revisionist
Integral Deferred Correction (RIDC).  RIDC methods, introduced in
\cite{explicit_ridc}, are ``parallel-across-the-step'' integrators
that can be efficiently implemented with multi-core
\cite{explicit_ridc,implicit_ridc}, multi-GPGPU \cite{ridc_gpgpu}, and
multi-node \cite{ridc_dd} architectures.  The ``revisionist''
terminology was adopted to highlight that (1) RIDC is a {revision} of
the standard integral defect correction (IDC) formulation~\cite{DGR00}, and (2)
successive corrections, running in parallel but (slightly) lagging in
time, {revise} and improve the approximation to the solution.

RIDC methods have been shown to be effective parallel time-integration
methods. They can typically produce a high-order solution in
essentially the same amount of wall-clock time as the constituent
lower-order methods. In general, for a given amount of wall-clock
time, RIDC methods are able to produce a more accurate solution than
conventional methods. These results have thus far been demonstrated
with constant time steps. It has long been accepted that local error
estimation and adaptive step-size control form a critical part of a
robust and efficient strategy for solving initial-value problems in
ordinary differential equations (ODEs), in particular problems with
multiple timescales; see e.g.,~\cite{ode1}. Accordingly, in order to
assess the practical viability of RIDC methods, it is important to
establish whether they can operate effectively with variable step
sizes. It turns out that there are subtleties associated with
modifying the RIDC framework to incorporate functionality for local
error estimation and adaptive step-size control: there are a number of
different implementation options, and some of them are more effective
than others.

The remainder of this paper is organized as follows.  In
Section~\ref{sec:Review}, we review the ideas behind RIDC as well as
strategies for local error estimation and step-size control.  We then
combine these ideas to propose various strategies for RIDC methods
with error and step-size control.  In Section~\ref{sec:adaptive_ridc},
we describe the implementation of these strategies within the RIDC
framework and suggest avenues that can be explored for a
production-level code.  In Section~\ref{sec:results}, we demonstrate
% that RIDC methods achieve their designed order of accuracy with
% non-uniform stepsizes and 
that the use of local error estimation and adaptive step-size control
inside RIDC is computationally advantageous.
%In particular, numerical
%experiments are presented that suggest that RIDC methods with error
%and step-size control are competitive with their standard Runge--Kutta
%(RK) counterparts. 
Finally, in Section~\ref{sec:conclusions}, we summarize the
conclusions reached from this investigation and comment on some
potential directions for future research.

\section{Review of Relevant Background}
\label{sec:Review}
We are interested in numerical solutions to initial-value problems
(IVPs) of the form
\begin{align}
  \label{eqn:ode}
  \left\{\begin{array} {l}
      \displaystyle
      y'(t)  = f(t,y(t)), \quad t \in[a, b], \\
      \displaystyle
      y(a) = y_a.
    \end{array}
  \right.
\end{align}
where $y(t) : \mathbb{R} \to \mathbb{R}^m$, $y_a\in\mathbb{R}^m$, and
$f:\mathbb{R}\times\mathbb{R}^m\to\mathbb{R}^m$. We first review RIDC
methods, a family of parallel time integrators that can be applied to
solve~\eqref{eqn:ode}.  Then, we review strategies for local error
estimation and adaptive step-size control for IVP solvers.

\subsection{RIDC}
RIDC methods \cite{explicit_ridc,implicit_ridc,ridc_gpgpu} are a class
of time integrators based on integral deferred correction \cite{DGR00}
that can be implemented in parallel via pipelining.  RIDC methods
first compute an initial (or {\em provisional}) solution, typically
using a standard low-order scheme, followed by one or more
corrections.  Each correction {revises} the current solution and
increases its formal order of accuracy.  After initial start-up costs,
the predictor and all the correctors can be executed in parallel.  It
has been shown that parallel RIDC methods give almost perfect 
parallel speedups~\cite{explicit_ridc}.  In this section, we review RIDC
algorithms, generalizing the overall framework slightly to allow for
non-uniform step-sizes on the different correction
levels.  % We shall assume presently that the time domain has
% already been divided into potentially non-uniform intervals on each
% level.

We denote the nodes for correction level $\ell$ by
\begin{align*}
 a = t_0^{[\ell]} < t_1^{[\ell]} < \cdots < t_{N^{[\ell]}}^{[\ell]} = b,
\end{align*}
where $N^{[\ell]}$ denotes the number of time steps on the level $\ell$.
In practice, the nodes on each level are obtained dynamically by the
step-size controller.

\subsubsection{The Predictor}
\label{sec:predictor}
To generate a provisional solution, a low-order integrator is applied
to solve the IVP~\eqref{eqn:ode}.  For example, a first-order forward
Euler integrator applied to~\eqref{eqn:ode} gives
\begin{align}
  \eta^{[0]}_{n} = \eta^{[0]}_{n-1} + \left(t_{n}^{[0]} - t_{n-1}^{[0]}\right)
  f(t_{n-1}^{[0]},\eta^{[0]}_{n-1}), 
  \label{eqn:pred}
\end{align}
for $n = 1,2,\ldots,N^{[0]}$, with $\eta_{0}^{[0]} = y_a$, and where we have
indexed the prediction level as level $0$.  We denote
$\eta^{[\ell]}(t)$ as a continuous extension \cite{ode1} of the
numerical solution at level $\ell$, i.e., a piecewise polynomial
$\eta^{[0]}(t)$ that satisfies
\begin{align*}
  \eta^{[0]}(t_n^{[0]}) = \eta_{n}^{[0]}.
\end{align*}
The continuous extension of a numerical solution is often of the same
order of accuracy as the underlying discrete solution \cite{ode1}.
Indeed, for the purposes
of this study, we assume $\eta^{[\ell]}(t)$ is of the same order as
$\eta^{[\ell]}_n$.

\subsubsection{The Correctors}
\label{sec:correctors}
Suppose an approximate solution $\eta(t)$ to IVP~\eqref{eqn:ode} is
computed.  Denote the exact solution by $y(t)$.  Then, the error of
the approximate solution is $ e(t) = y(t)-\eta(t)$.  If we define the
defect as $\delta(t) = f(t,\eta(t))- \eta'(t)$, then
\begin{align*}
  e'(t) = y'(t) - \eta'(t) = f(t,\eta(t)+e(t)) - f(t,\eta(t)) +
  \delta(t).
\end{align*}
The error equation can be written in the form
\begin{align}
  \label{eqn:error_eqn}
  \left[e(t)-\int_a^t\delta(\tau)\,d\tau\right]' =
  f\left(t,\eta(t)+e(t)\right) - f\left(t,\eta(t)\right)
\end{align}
subject to the initial condition $e(a) = 0$. In RIDC, the corrector at
level $\ell$ solves for the error $e^{[\ell-1]}(t)$ of the solution
$\eta^{[\ell-1]}(t)$ at the previous level to generate the corrected
solution $\eta^{[\ell]}(t)$,
\begin{align*}
  \eta^{[\ell]}(t) = \eta^{[\ell-1]}(t) + e^{[\ell-1]}(t).
\end{align*}
For example, a corrector at level $\ell$ that corrects
$\eta^{[\ell-1]}(t)$ by applying a first-order forward Euler
integrator to the error equation~\eqref{eqn:error_eqn} takes the form
\begin{multline*}
  e^{[\ell-1]}(t_{n}^{[\ell]}) -e^{[\ell-1]}(t_{n-1}^{[\ell]}) -
  \int_{t_{n-1}^{[\ell]}}^{t_n^{[\ell]}} 
  \delta^{[\ell-1]}(\tau)\,d\tau = \\  
 \Delta t_{n}^{[\ell]}
  \left[f\left(t_{n-1}^{[\ell]},\eta^{[\ell-1]}(t_{n-1}^{[\ell]}) +
    e^{[\ell-1]}(t_{n-1}^{[\ell]}) \right)
    - f\left(t_{n-1}^{[\ell]},\eta^{[\ell-1]}(t_{n-1}^{[\ell]})  \right)
 \right],
\end{multline*}
where $ \Delta t_{n}^{[\ell]} = t_{n}^{[\ell]} - t_{n-1}^{[\ell]}$.
After some algebraic manipulation, one obtains
\begin{align}
  \label{eqn:corr}
  \eta^{[\ell]}_n = \,& \eta^{[\ell]}_{n-1} + \Delta t_{n}^{[\ell]}\left(
  f\left(t_{n-1}^{[\ell]},\eta^{[\ell]}(t_{n-1}^{[\ell]})  \right)
    - f\left(t_{n-1}^{[\ell]},\eta^{[\ell-1]}(t_{n-1}^{[\ell]})  \right)
  \right)\\
  \nonumber
  &+   \int_{t_{n-1}^{[\ell]}}^{t_n^{[\ell]}} 
  f\left(\tau,\eta^{[\ell-1]}(\tau)\right)\,d\tau.
\end{align}
The integral in equation~\eqref{eqn:corr} is approximated using
quadrature,
\begin{align}
  \label{eqn:quadrature}
  \int_{t_{n-1}^{[\ell]}}^{t_{n}^{[\ell]}}
  f\left(\tau,\eta^{[\ell-1]}(\tau)\right)\,d\tau
  \approx
  \sum_{i=1}^{|\vec{\mathcal{T}}^{[\ell]}_n|} \alpha_{n,i}^{[\ell-1]}
  f\left(\tau_i,\eta^{[\ell-1]}(\tau_i)\right), 
  \quad \tau_i \in \vec{\mathcal{T}}^{[\ell]}_n,
\end{align}
where the set of quadrature nodes, $\vec{\mathcal{T}}^{[\ell]}_n$,
for a first-order corrector satisfies
\begin{enumerate}
\item $|\vec{\mathcal{T}}^{[\ell]}_n| = \ell+1 $
\item $\vec{\mathcal{T}}^{[\ell]}_n \subseteq \{t_n^{[\ell-1]} \}_{n=0}^{N^{[\ell-1]}}$
\item $\min (\vec{\mathcal{T}}^{[\ell]}_n)  \le t_{n-1}^{[\ell]}$
\item $\max (\vec{\mathcal{T}}^{[\ell]}_n)  \ge t_n^{[\ell]}$
\end{enumerate}
The quadrature weights, $\alpha_{n,i}^{[\ell-1]}$, are found by
integrating the interpolating Lagrange polynomials exactly,
\begin{align}
  \alpha_{n,i}^{[\ell-1]} = \prod_{j=1,j\neq i} ^{|\vec{\mathcal{T}}^{[\ell]}_n|}
  \int_{t_{n-1}^{[\ell]}}^{t_n^{[\ell]}}
  \frac{(t-\tau_j)}{(\tau_i - \tau_j)}\,dt,
    \quad \tau_i \in \vec{\mathcal{T}}^{[\ell]}_n.
    \label{eqn:quadrature_weights}
\end{align}
The term $f\left(t_{n-1}^{[\ell]},\eta^{[\ell-1]}(t_{n-1}^{[\ell]})\right)$ in
equation~\eqref{eqn:corr} is approximated using Lagrange
interpolation,
\begin{align}
  f\left(t_{n-1}^{[\ell]},\eta^{[\ell-1]}(t_{n-1}^{[\ell]})\right)
  \approx
  \sum_{i=1}^{|\vec{\mathcal{T}}^{[\ell]}_n|} \gamma_{n,i}^{[\ell-1]}
  f\left(\tau_i,\eta^{[\ell-1]}(\tau_i)\right), 
  \quad \tau_i \in \vec{\mathcal{T}}^{[\ell]}_n,
  \label{eqn:interpolation}
\end{align}
where the same set of nodes, $\vec{\mathcal{T}}^{[\ell]}_n$, for the
quadrature is used for the interpolation.  The interpolation weights
are given by
\begin{align}
  \gamma_{n,i}^{[\ell-1]} = \prod_{j=1,j\neq i} ^{|\vec{\mathcal{T}}^{[\ell]}_n|}
  \frac{(t_{n-1}^{[\ell]}-\tau_j)}{(\tau_i - \tau_j)},
  \quad \tau_i \in \vec{\mathcal{T}}^{[\ell]}_n.
  \label{eqn:interpolation_weights}
\end{align}

%.  The number of nodes needed to accurately approximate the
%integral is related

%For the $p$\textsuperscript{th} correction loop, $p$ nodes
%are needed in the stencil to accurately approximate the integral.
%There are various choices for the stencil, but at the very minimum,
%the stencil chosen from the nodes $\{t_i^{[p-1]}\}$ should
%enclose the nodes $t_n^{p]}$ and $t_{n+1}^{[\ell]}$.  Similarly, the term
%  $ f\left(t_{n-1}^{[\ell]},\eta^{[p-1]}(t_{n-1}^{[\ell]})\right)$ is
%  approximated by finding a stencil containing $p$ nodes from
%  $\{t_i^{[p-1]}\}$  that contains $t_n^{[\ell]}$, and using
%    that stencil to form and evaluate an interpolating polynomial at
%    the node $t_n^{[\ell]}$.

\subsection{Adaptive Step-Size Control}
Adaptive step-size control is typically used to achieve a
user-specified error tolerance with minimal computational effort by
varying the step-sizes used by an IVP integrator.  This is commonly
done based on a local error estimate.  It may also be desirable that
the step-size vary smoothly over the course of the integration.  We
review common techniques for estimating the local error, followed by
algorithms for optimal step-size selection.

\subsubsection{Error Estimators}
\label{sec:error_estimate}
Two common approaches for estimating the local truncation error of a
single-step IVP solver are through the use of Richardson extrapolation
(commonly used within a step-size selection framework known as
step-doubling) and embedded Runge--Kutta pairs \cite{ode1}.
Step-doubling is perhaps the more intuitive technique.  The solution
after each step is estimated twice: once as a full step and once as
two half steps.  The difference between the two numerical estimates
gives an estimate of the truncation error.  For example, denoting the
exact solution to IVP~\eqref{eqn:ode} at time $t_n+\Delta t$ as
$y(t_n+\Delta t)$, the forward Euler step starting from the exact
solution at time $t_n$ and using a step size of size $\Delta t$ is
\begin{align*}
  \eta_{1,n+1} = y(t_n) + \Delta t \,f(t_n,y_n),
\end{align*}
and the forward Euler step using two steps of size
$\frac{\Delta t}{2}$ is
\begin{align*}
  \eta_{2,n+1} = \left(y(t_n) + \frac{\Delta t}{2} \,f(t_n,y_n)\right) +
  \frac{\Delta t}{2} f\left(t_n+\frac{\Delta t}{2},y(t_n) +
  \frac{\Delta t}{2} \,f(t_n,y_n) \right).
\end{align*}
Because forward Euler is a first-order method (and thus has a local
truncation error of $\mathcal{O}(\Delta t^2)$), the two numerical
approximations satisfy
\begin{align*}
  y(t_n+\Delta t) &= \eta_{1,n+1} + (\Delta t)^2\phi 
  + \mathcal{O}(\Delta t^3) + \cdots, \\
  y(t_n+\Delta t) &= \eta_{2,n+1} + 2\left(\frac{\Delta t}{2}\right)^2\phi 
  + \mathcal{O}(\Delta t^3) + \cdots, 
\end{align*}
where a Taylor series expansion gives that $\phi$ is a constant
proportional to $y''(t_n)$.  The difference between the two numerical
approximations gives an estimate for the local truncation error of
$\eta_{2,n+1}$,
\begin{align*}
  e_{n+1} = \eta_{2,n+1} - \eta_{1,n+1} = \frac{\Delta t^2}{2}\phi +
  \mathcal{O}(\Delta t^3).
\end{align*}

An alternative approach to estimating the local truncation error is to
use embedded RK pairs \cite{embedded_rk}.  An $s$-stage Runge--Kutta method is
a single-step method that takes the form
\begin{align*}
  \eta_{n+1} &= \eta_n + \Delta t \sum_{i=1}^s b_i k_i,
  \\
  \hspace*{-15mm}\text{where}
  \qquad \qquad \
  k_i &= f\left(t_i+c_ih,\eta_n + \Delta t \sum_{j=1}^s a_{ij} k_j
  \right), \quad i=1,2,\ldots,s.
\end{align*}
The idea is to find two single-step RK methods, typically one with
order $p$ and the other with order $p-1$, that share most (if not all)
of their stages but have different quadrature weights.  This is represented
compactly in the extended Butcher tableau
  \begin{equation*}
    \begin{array}{l|c}
      c & A \\
      \hline
      & b\\
      & \hat{b}
    \end{array}
  \end{equation*}
  Denoting the solution from the order-$p$ method as
\begin{subequations}
\begin{equation}
\label{eqn:eta*}
\eta^*_{n+1} = \eta_n + \Delta t \sum_{i=1}^s \hat{b}_i k_i,
\end{equation}
and the solution from the order-$(p-1)$ method as
\begin{equation}
\label{eqn:eta}
\eta_{n+1} = \eta_n + \Delta t \sum_{i=1}^s b_i k_i,
\end{equation}
the error estimate is
\begin{align}
  e_{n+1} =\eta_{n+1} - \eta^*_{n+1} = \Delta t \sum_{i=1}^s (b_i-\hat{b}_i) k_i,
  \label{eqn:double_truncation}
\end{align}
\end{subequations}
which is $\mathcal{O}(\Delta t^p)$.

A third approach for approximating the local truncation error is
possible within the deferred correction framework.  We observe that
that in solving the error equation~\eqref{eqn:error_eqn}, one is in
fact obtaining an approximation to the error.  As discussed in
Section~\ref{sec:using_error_equations}, it can be shown that the approximate
error after $\ell$ first-order corrections satisfies $o(\Delta
t^{p_0+\ell+1})$.  We shall see in Section~\ref{sec:using_error_equations}
that this error estimate proves to be a poor choice for optimal step
size selection because in our formulation the time step selection for
level $\ell$ does not allow for the refinement of time steps at
earlier levels.

% Because~\eqref{eqn:error_eqn} is
%formulated from~\eqref{eqn:ode}, the error estimate obtained is a
%global error measure.  This contrasts with the local truncation error
%estimates from step-doubling or embedded RK pairs.  \todo{Note: this
%  seems like a big deal, make sure we follow-up on this point later.}
%\ongbw{This paragraph is not worded correctly. I think the error
%  estimate is still local, however, the convergence problem arises
%  because we're not taking time steps uniformly to zero on all levels}

\subsubsection{Optimal Step-Size Selection}
Given an error estimate from Section~\ref{sec:error_estimate} for a
step $\dt$, one would like to either accept or reject the step based
on the error estimate and then estimate an optimal step-size for the
next time step or retry the current step.  Following~\cite{ode2},
Algorithm~\ref{alg:step_size_selection} outlines optimal step-size
selection given an estimate of the local truncation
error.  In lines 1--4, one computes a scaled error estimate.
% using the
%differences between the regular-step and the double-step solutions.
In line 5, an optimal time step is computed by scaling the current
time step.  In lines 6--10, a new time step is suggested; a more
conservative step-size is suggested if the previous step was rejected.
%\todo{Not sure we need this algorithm repeated in the  paper.}
\begin{algorithm}[htbp]
  \DontPrintSemicolon
  \SetAlgoLined

  \KwIn{\\ $y_n$: approximate solution at time $t_n$; \\
    $y_{n+1}$: approximate solution at time $t_{n+1}$; \\
    $e_{n+1}$: error estimate for $y_{n+1}$; \\
    $p$: order of integrator; \\
    $m$: number of ODEs; \\
    {\tt atol}, {\tt rtol}: user specified tolerances; \\
    {\tt prev\_rej}: flag that indicates whether the previous step was rejected; \\
    $\alpha < 1$: safety factor; \\
    $\beta > 1$: allowable change in step-size.}

  \KwOut{\\{\tt accept\_flag}: flag to accept or reject this step; \\$\dt_{new}$:
    optimal time step}

  % Compute error estimate, $e(i) = |\eta_2(i)-\eta_1(i)|/(2^p-1),\, i=1,2,\ldots,m$.

  Set $a(i) = \max\{|y_n(i)|,|y_{n+1}(i)|\},\, i=1,2,\ldots,m$.

  Compute $\tau(i) = {\tt atol} + {\tt rtol}*a(i),\, i=1,2,\ldots,m$.

  Compute $\epsilon = \sqrt{\frac{\sum_{i=1}^m (e(i)/\tau(i))^2}{m}}$.

  Compute $\dt_{opt} = \dt\,(\frac1\epsilon)^{1/(p+1)}$.

  \eIf{{\tt prev\_rej}} { 
    $\dt_{new} = \alpha\min\{\dt,\max\{\dt_{opt},\dt/\beta\}\}$
  } {
    $\dt_{new} = \alpha\min\{\beta\dt,\max\{\dt_{opt},\dt/\beta\}\}$
  }

  \eIf{$\epsilon > 1$} {
   $ {\tt accept\_flag} = 1$
   } { 
    ${\tt accept\_flag} = 0$
  }
  \caption{Optimal step-size selection algorithm. The approximate
    solution, the error estimate, and its order are provided as
    inputs.  For the numerical experiments in Section~\ref{sec:results},
    we fix $\alpha=0.9$, $\beta=10$.}
  \label{alg:step_size_selection}
\end{algorithm}

\section{RIDC with Adaptive Step-Size Control}
\label{sec:adaptive_ridc}
There are numerous adaptive step-size control strategies that can be
implemented within the RIDC framework.  We consider three of them in
this paper as well as discuss other strategies that are possible.

\subsection{Adaptive Step-Size Control: Prediction Level Only}
\begin{figure}[htbp]
  \centering
  \includegraphics{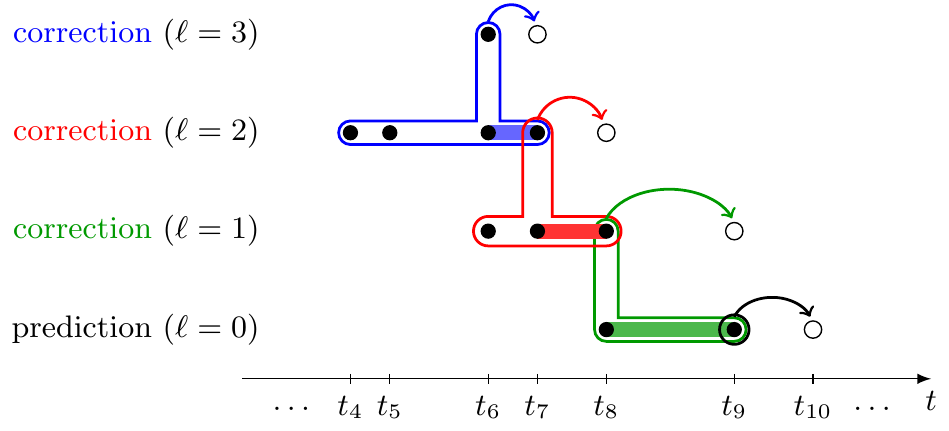}
  \caption{Schematic diagram of step-size control on the prediction
    level only.  The filled circles denote previously computed and
    stored solution values at particular times.  The corrections are
    run in parallel (but lagging in time) and the open circles
    indicate which values are being simultaneously computed.
    The stencil of points required by each level is shown by the
    ``bubbles'' surrounding certain grid points; the thick
    horizontal shading indicates the integrals needed in
    \eqref{eqn:corr}.
    Note each level uses the same grid in time.
  }
  \label{fig:pred_only}
\end{figure}

One simple approach to step-size control with RIDC is to perform
adaptive step-size control on the prediction level only, for example
using step-doubling or embedded RK pairs as error estimators for the step-size control
strategy.  The subsequent correctors then use this grid unchanged
(i.e., without performing further step-size control).  With this strategy,
corrector $\ell$ is lagged behind corrector $\ell-1$ so that each node
simultaneously computes an update on its level (after an initial
startup period).  This is illustrated graphically in
\autoref{fig:pred_only}.  In principle, parallel speedup is maintained
with this approach, and a minimal memory footprint is required in an
implementation.  Additionally, an interpolation step is circumvented
because the nodes are the same on each level.  There are however a few
potential drawbacks in this approach. First, it is not clear how to
distribute the user-defined tolerance among the levels. Clearly,
satisfying the user-specified tolerance on the prediction level
defeats the purpose of the deferred correction approach.
Estimating a reduced tolerance criterion may be possible a priori, but
such an estimate would at present be ad hoc.  Second, there is no
reason to expect the corrector~\eqref{eqn:corr} should take the same
steps to satisfy an error tolerance when computing a numerical
approximation to the error equation~\eqref{eqn:error_eqn}.

\subsection{Adaptive Step-Size Control: All Levels}
\label{sec:adapt_all}
\begin{figure}[htbp]
  \centering
  \includegraphics{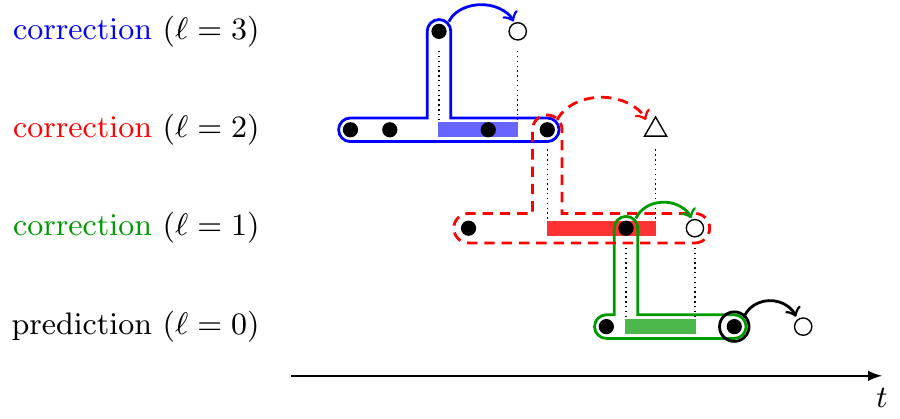}
  \caption{Schematic diagram of a scenario when step-size control is
    applied on all levels.  Unlike in Figure~\ref{fig:pred_only}, here
    each level has its own grid in time.  Solid circles indicate
    particular times and levels where the solution is known.  In this
    particular diagram, levels $\ell=0,1,3$ are all able to advance
    simultaneously to the open circles.  However, correction level
    $\ell=2$ is unable to advance to the time indicated by the
    triangle symbol because correction level $\ell=1$ has not yet
    computed far enough.  The stencil of points required by each level
    is shown by the ``bubbles'' surrounding certain grid points; the
    thick horizontal shading indicates the integrals needed in
    \eqref{eqn:corr}.  Note in particular that the dashed stencil
    includes a open circle at level $\ell=1$ that is not yet
    computed.  }
  \label{fig:all_levels}
\end{figure}

A generalization of the above formulation is to utilize adaptive
step-size control to solve the error equations~\eqref{eqn:error_eqn}
as well.  The variant we consider is step-doubling on all levels,
where each predictor and corrector performs
Algorithm~\ref{alg:step_size_selection}; embedded RK pairs can also be
used to estimate the error for step-size adaptivity on all levels.
Intuitively, step-size control on every level gives more opportunity
to detect and adapt to error than simply adapting using the
(lowest-order) predictor.  For example, this allows the corrector take
a smaller step if necessary to satisfy an error tolerance when solving
the error equation.  Some drawbacks are: (i) an interpolation step is
necessary because the nodes are generally no longer in the same
locations on each level, (ii) more memory registers are required, and
(iii) there is a potential loss of parallel efficiency because a
corrector may be stalled waiting for an adequate stencil to become
available to compute a quadrature approximation to the integral in
equation~\eqref{eqn:corr}.  Another issue --- both a potential benefit
and a potential drawback --- is the number of parameters that can be
tuned for each problem.  A discussion on the effect of tolerance
choices for each level is provided in Section~\ref{sec:results}.  One
can in practice also tune step-size control parameters $\alpha,
\beta$, {\tt atol}, and {\tt rtol} for
Algorithm~\ref{alg:step_size_selection} separately on each level.
Figure~\ref{fig:all_levels} highlights that some nodes might not be
able to compute an updated solution on their current level if an
adequate stencil is not available to approximate the integral in
equation~\eqref{eqn:corr} using quadrature.  In this example, the
level $\ell=2$ correction is unable to proceed, whereas the prediction
level $\ell=0$ and corrections $\ell=1$ and $\ell=3$ are all able to
advance the solution by one step.

%\subsection{Embedded RK Pairs at all Levels}
%
%The algorithm that uses embedded RK pairs for error and step-size
%control proceeds as Algorithm~\ref{alg:step_size_selection} with the
%modification that $u_1$ and $u_2$ are taken to be $\eta_{n+1}^*$ and
%$\eta_{n+1}$ from (\ref{eqn:eta*}) and (\ref{eqn:eta}), respectively,
%and the error estimate is computed from (\ref{eqn:double_truncation}).
%Although only a first-order accuracy improvement can be expected by
%using any embedded RK pair (which necessarily require at least two
%stages), there might be stability advantages over using a step doubled
%solution which might make such a formulation attractive.  
%
%\todo{C\&R: yes we like this, Ray will take a look.  We definitely
%think the stability comment makes sense.}

\subsection{Adaptive Step-Size Control: Using the Error Equation}
\label{sec:using_error_equations}

A third strategy one might consider is adaptive step-size control for
the error equation~\eqref{eqn:error_eqn} using the solution to the
error equation itself as the error estimate.  (One still uses
step-doubling or embedded RK pairs to obtain an error estimate for step-size
control on the predictor equation~\eqref{eqn:ode}). At first glance,
this looks promising provided the order of the integrator can be
established because it is used to determine an optimal step-size. One
would expect computational savings from utilizing available error
information, as opposed to estimating it via step-doubling or an
embedded RK pair.

If first-order predictor and first-order correctors are used to
construct the RIDC method, the analysis in \cite{XXS} can be easily extended
to the proposed RIDC methods with adaptive step-size
control. We note that the numerical quadrature approximation given in
equation~\eqref{eqn:quadrature} and the numerical interpolation given
in equation~\eqref{eqn:interpolation} are accurate to the order
$\mathcal{O}(\Dt_n^{\ell+2})$; this is sufficient for the inductive proof
in \cite{XXS} to hold.  Hence, one can show that the method has a
formal order of accuracy $\mathcal{O}(\Dt^{\ell+2})$, where $\Dt =
\max_{n,\ell}(t_{n}^{[\ell]}-t_{n-1}^{[\ell]})$.

Although the formal order of accuracy can be established, using the
error estimate from successive levels is a poor choice for optimal
step-size selection.  Consider step-size selection for level $\ell$,
time step $t_n^{[\ell]}$, using $\eta_n^{[\ell]} -
\eta^{[\ell-1]}(t_n^{[\ell]})$ as the error estimator in Algorithm 1.
The optimal step-size is chosen to control the local error estimate
via the step-size $\Dt_n^{[\ell]} = t_{n}^{[\ell]}-t_{n-1}^{[\ell]}$.
However, the local error for the correctors generally contains
contributions from the solutions at all the previous levels. The
validity of the asymptotic local error expansion of the RIDC method in
terms of $\Dt_n^{[\ell]}$ requires that $\Dt =
\max_{n,\ell}(t_{n}^{[\ell]}-t_{n-1}^{[\ell]})$ be sufficiently small,
and it is not normally possible to guarantee this in the context of an
IVP solver. In other words, the step-size controller for a corrector
at a given level cannot control the entire local error, and hence
standard step-control strategies, which are predicated on the validity
of error expansions in terms of only the step-size to be taken, cannot
be expected to perform well.  We present some numerical tests in
Section~\ref{sec:numerics_successive} to illustrate the difficulties
with using successive errors as the basis for step-size control.

%% Whilst this might seem like a reasonable
%% idea at first glance, some analysis reveals that
%% \begin{align*}
%%   e^{[\ell]}(t_{n+1}^{[\ell]}) -  e^{[\ell]}(t_{n}^{[\ell]}) 
%%   = \mathcal{O}((t_{n+1}^{[\ell]}-t_{n}^{[\ell]})^p) + 
%%   \mathcal{O}((\Delta t^{[p-1]})^{\ell-1}),
%% \end{align*}
%% \todo{Need to define $\Delta t^{[p-1]}$. Why not use $\Delta
%%   t^{[\ell]}_{n+1}$?} \ongbw{see previous comment near line 94}
%% \todo{Ray!}
%% that
%% is, $e^{[\ell]}(t_{n+1}^{[\ell]}) - e^{[\ell]}(t_{n}^{[\ell]}) \to
%% \mathcal{O}((\Delta t^{[\ell-1]})^{p-1})$ as $t_{n+1}^{[\ell]}\ \to
%% t_{n}^{[\ell]}$ because interpolation error from using solutions from
%% previously computed levels necessarily exists
%% in~\eqref{eqn:error_eqn}.  Although this interpolation error also
%% exists in the formulation that uses step-doubling at every level,
%% Section~\ref{sec:adapt_all}, it is more problematic in this
%% formulation.  The reason is quite subtle: when using step-doubling to
%% approximate the error for the error equation~\eqref{eqn:error_eqn},
%% both numerical approximations ``contain'' the interpolation error, so
%% the difference between the two solutions still converges to 0 as $\dt
%% \to 0$.  This contrasts to using the solution to the error equation as
%% the estimate, where $e^{[\ell]}(t_{n+1}^{[\ell]}) -
%% e^{[\ell]}(t_{n}^{[\ell]}) \to \mathcal{O}((\Delta t^{[p-1]})^{p-1})$
%% as $t_{n+1}^{[\ell]}\ \to t_{n}^{[\ell]}$.  \todo{I'm not sure I get
%%   it, need to go over this carefully. It's probably important enough
%%   to write out.}

\subsection{Further Discussion}
There are many other strategies/implementation choices that affect the
overall performance of the adaptive RIDC algorithm.  Some have already
been discussed in the previous section.  We summarize
some of the implementation choices that must be made:
\begin{itemize}
\item The choice of how to estimate the error of the discretization
  must be made.  Three possibilities have already been mentioned:
  step-doubling, embedded RK pairs, and solutions to the error
  equation~\eqref{eqn:error_eqn}.  A combination of all three is also
  possible.
\item If an IVP method with adaptive step-size control is used to solve
  equation~\eqref{eqn:error_eqn}, choices must be made as to how the
  tolerances and step-size control parameters, $\alpha$ and $\beta$,
  are to be chosen for each correction level.

\end{itemize}

We also list a few implementation details that should be considered
when designing adaptive RIDC schemes.
\begin{itemize}
\item If adaptive step-size control is implemented on all levels, some
  correction levels may sit idle because the information required to
  perform the quadrature and interpolation in \eqref{eqn:corr} is not
  available.  This idle time adversely affects the parallel efficiency
  of the algorithm.  One possibility to decrease this idle time is
  instead of taking an ``optimal step'' (as suggested by the step-size
  control routine), one could take a smaller step for which the
  quadrature and interpolation stencil is available.  There is some
  flexibility in choosing exactly which points are used in the
  quadrature stencil, and it might also be possible to choose a stencil
  to minimize the time that correction levels are sitting idle.
% \todo{Colin changed this item and added the one below.  Does this
%     address the points we're trying to make here? -- yes (ongbw)}

\item Because values are needed from lower-order correction levels,
  the storage required by a RIDC scheme depends on when values can be
  overwritten (see, e.g., the stencils in Figures~\ref{fig:pred_only}
  and~\ref{fig:all_levels}).  Thus to avoid increasing the storage
  requirements, the prediction level and each correction level should
  not be allowed to get too far ahead of higher correction levels.
  Although this is also the case for the non-adaptive RIDC schemes
  \cite{explicit_ridc,implicit_ridc}, if adaptive step-size control is
  implemented on all levels (\autoref{fig:all_levels}), the memory
  footprint is likely to increase.  Some consideration should thus be
  given to a potential trade-off between parallel efficiency and the
  overall memory footprint of the scheme.

\item It is important to reduce round-off error when computing the
  quadrature weights~\eqref{eqn:quadrature_weights} and the
  interpolation weights~\eqref{eqn:interpolation_weights}.  This
  can be done by through careful scaling and control of the order of
  the floating-point operations \cite{bradie2006friendly}.

\item If one wishes to use higher-order correctors and predictors to
  construct RIDC integrators, we note that the convergence analysis in
  \cite{IDC2,IDC1,idc-ark} only holds for uniform steps.  A
  non-uniform mesh introduces discrete ``roughness'' (see
  \cite{IDC1}); hence, an increase of only one order per correction
  level is guaranteed even though a high-order method is used to solve
  equation~\eqref{eqn:error_eqn}. 
  %Otherwise, the discussion in
  %Section~\ref{sec:adapt_all} on the potential benefits or drawbacks
  %of adaptive step-size control on all levels applies to this third
  %strategy. Additionally, 

\end{itemize}

Additionally, the RIDC framework, by construction, solves a series of
error equations to generate a successively more accurate solution.
This framework can be potentially be exploited to generate {\em
  order-adaptive} RIDC methods.  For example, one might control the
number of corrector levels adaptively based on an error estimate.

\section{Numerical Examples}
\label{sec:results}
We focus on the solutions to three nonlinear IVPs.  The first is
presented in \cite{AHKW}; we refer to it as the Auzinger IVP:
\begin{align}
  \tag{AUZ}
  \label{eqn:auz}
  \left\{\begin{array} {l}
      \displaystyle
      \vspace*{1mm}
      y_1' = -y_2 + y_1(1-y_1^2 - y_2^2),\\
      \vspace*{1mm}
      \displaystyle
      y_2' = y_1 + 3y_2(1-y_1^2-y_2^2),\\
      \displaystyle
      y(0) = (1,0)^T, \quad t \in[0, 10], \\
    \end{array}
  \right.
\end{align}
that has the analytic solution $y(t) = (\cos{t},\sin{t})^T$. 

The second is the IVP associated with the Lorenz attractor:
\begin{align}
  \tag{LORENZ}
  \label{eqn:lorenz}
  \left\{\begin{array} {l}
      \displaystyle
      \vspace*{1mm}
      y_1' = \sigma(y_2-y_1),\\
      y_2' = \rho y_1 - y_2 - y_1y_3, \\
      \vspace*{1mm}
      y_3' = y_1y_2 - \beta y_3,\\
      y(0) = (1,1,1)^T, \quad t \in[0, 1]. 
    \end{array}
  \right.
\end{align} 
For the parameter settings $\sigma = 10,\ \rho = 28,\ \beta = 8/3$,
this system is highly sensitive to perturbations, and an IVP
integrator with adaptive step-size control may be advantageous.

The third is the restricted three-body problem from~\cite{ode1};  we refer
to it as the Orbit IVP:
\begin{align}
  \tag{ORBIT}
  \label{eqn:orbit}
  \left\{\begin{array} {l}
      \displaystyle
      \vspace*{1mm}
      \displaystyle
      y_1'' = y_1 + 2y_2' - \mu'\frac{y_1+\mu}{D_1} - \mu\frac{y_1-\mu'}{D_2}\\
      \vspace*{1mm}
      \displaystyle
      y_2'' = y_2 - 2y_1' - \mu'\frac{y_2}{D_1} - \mu\frac{y_2}{D_2}\\
      \vspace*{1mm}
      \displaystyle
      D_1 = \left((y_1+\mu)^2+y_2^2\right)^{3/2}, \quad
      D_2 = \left((y_1-\mu')^2+y_2^2\right)^{3/2}, \\
      \displaystyle
      \mu = 0.012277471, \quad \mu' = 1-\mu.
    \end{array}
  \right.
\end{align}
Choosing the initial conditions
\begin{align*}
  y_1(0) = 0.994, \quad y_1'(0) = 0, \quad y_2(0) = 0,\\
  y_2'(0) = -2.00158510637908252240537862224,
\end{align*}
gives a periodic solution with period $ t_{end} =
17.065216560159625588917206249$.

We now present numerical evidence to demonstrate that:
\begin{enumerate}
\item RIDC integrators with non-uniform step-sizes converge and
  achieve their designed orders of accuracy.
\item RIDC methods with adaptive step-size based on step-doubling and
  embedded RK error estimators, on the prediction level only, converge.
\item RIDC methods with adaptive step-size control based on
  step-doubling to estimate the local error on the prediction and
  correction levels converge; however, the step-sizes selected are poor
  (many rejected steps), even for the smooth Auzinger problem.
\item RIDC methods with adaptive step-size control based on
  step-doubling to estimate the local error on the prediction level
  but using the solution to the error equation for step-size control
  results is problematic.
\end{enumerate}

\subsection{RIDC with non-uniform step-sizes}
\begin{figure}[htbp]
  \centering
  \subfloat[Uniform steps]{
    \includegraphics{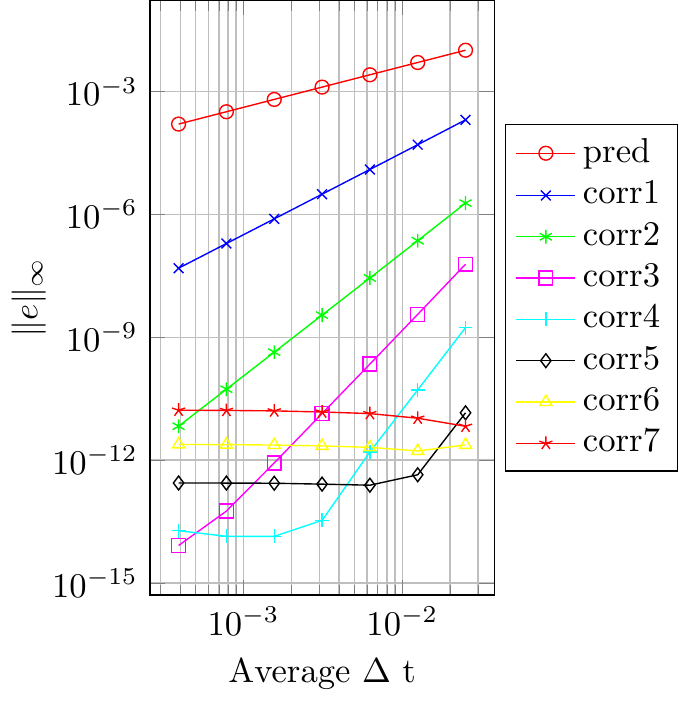}
  }
  \subfloat[Random steps, $\omega = 2$]{
    \includegraphics{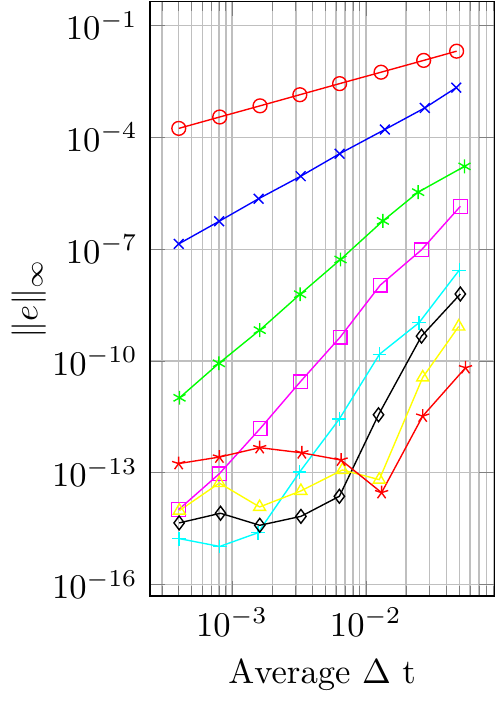}
  }\\
  \subfloat[Random steps, $\omega = 4$]{
    \includegraphics{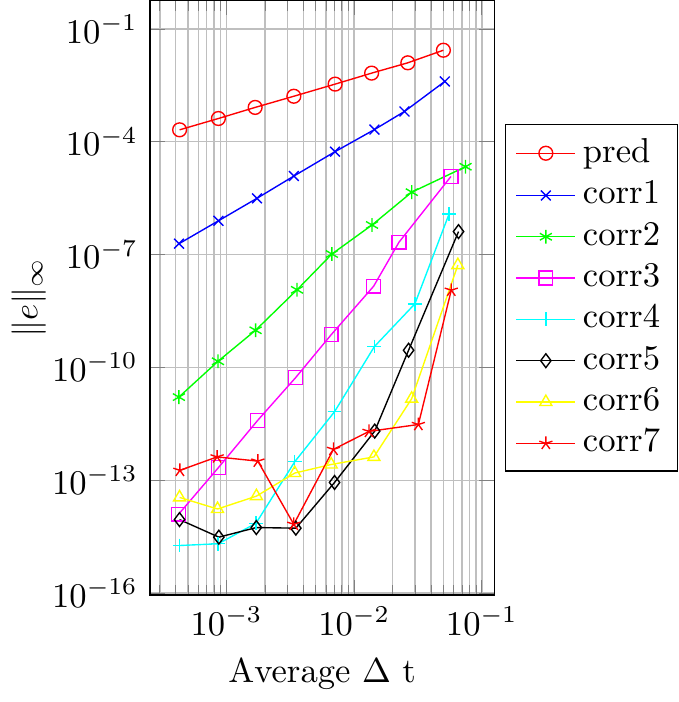}}
    \subfloat[Random steps, $\omega= 100$]{
    \includegraphics{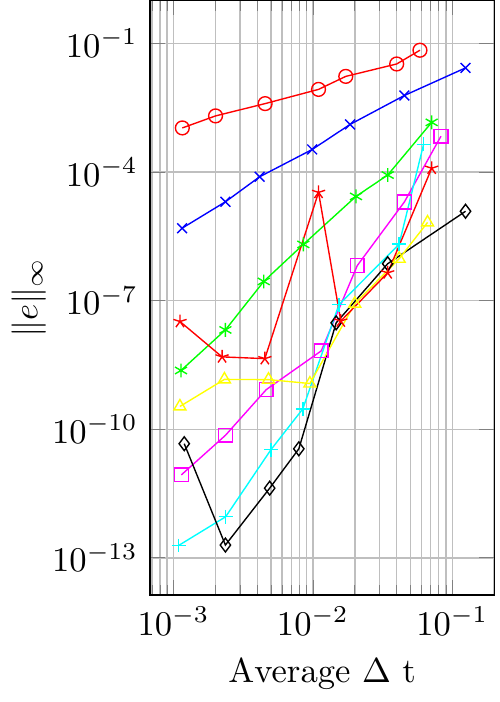}
  }
  \caption{Auzinger IVP: The design order is illustrated for the RIDC
    methods.  }
    \label{fig:auz}
\end{figure}

For our first numerical experiment, we demonstrate that RIDC
integrators with non-uniform step-sizes converge and achieve their
design orders of accuracy.  
% The grid is specified randomly as
%follows.  For a given $\Delta t^*$, set
%\begin{align*}
%  t_0^{[\ell]} = 0, \quad \Delta t_{-1}^{[\ell]}= \Delta t^*, \quad p=0,\ldots,P.
%\end{align*}
%Then, generate the non uniform grid on each level using
%\begin{align*}
%  t_{n+1}^{[\ell]} = t_n^{[\ell]} + \Delta t_{n}^{[\ell]}, \quad p=0,\ldots,P,
%  \quad n = 0, 1, \ldots
%\end{align*}
%where $\Delta t_n^{[\ell]}$ is randomly selected from $\left[\frac1\omega
%  \Delta t_{n-1}^{[\ell]}, \omega\Delta t_{n-1}^{[\ell]} \right]$ with
%$\omega \ge 1$.  Here, $\omega$ controls how rapidly steps are allowed
%to increase or decrease.  One chooses the last step for each level
%such that the domain of integration is exactly $[0,T]$.
%
Figure~\ref{fig:auz} shows the classical convergence study (error
as a function of mean step-size) for the RIDC integrator applied to
equation~\eqref{eqn:auz}.  Figure~\ref{fig:auz}(a) shows the
convergence of RIDC integrators with uniform step-sizes;
Figures~\ref{fig:auz}(b)--(d) show the convergence of RIDC integrators
when random step-sizes are chosen.  The random step-sizes are chosen so
that
\begin{align*}
  \dt_n^{[\ell]} \in \left[\frac1\omega \Delta t_{n-1}^{[\ell]},\omega\Delta t_{n-1}^{[\ell]} \right],
  \quad \omega\in \mathbb{R},
\end{align*}
where $\omega$ controls how rapidly a step-size is allowed to change.
The figures show that RIDC integrators with non-uniform step-sizes
achieve their designed order of accuracy (each additional correction
improves the order of accuracy by one), at least up to order 6.  In
\autoref{fig:auz} (corresponding to RIDC with uniform step-sizes),
we observe that the error stagnates at a value significantly larger
than machine precision.  This is likely due to numerical issues
associated with quadrature on equispaced nodes
\cite{rational_integration}.  We note that $\omega=1$ gives the
uniformly distributed case.  We also observe that as the ratio of the
largest to the smallest cell increases, the performance of
higher-order RIDC methods degrades, likely due to round-off error
associated with calculating the quadrature and interpolation weights.

\begin{figure}[htbp]
  \centering
  \subfloat[Ratio = 1 (uniform)]{
  \includegraphics{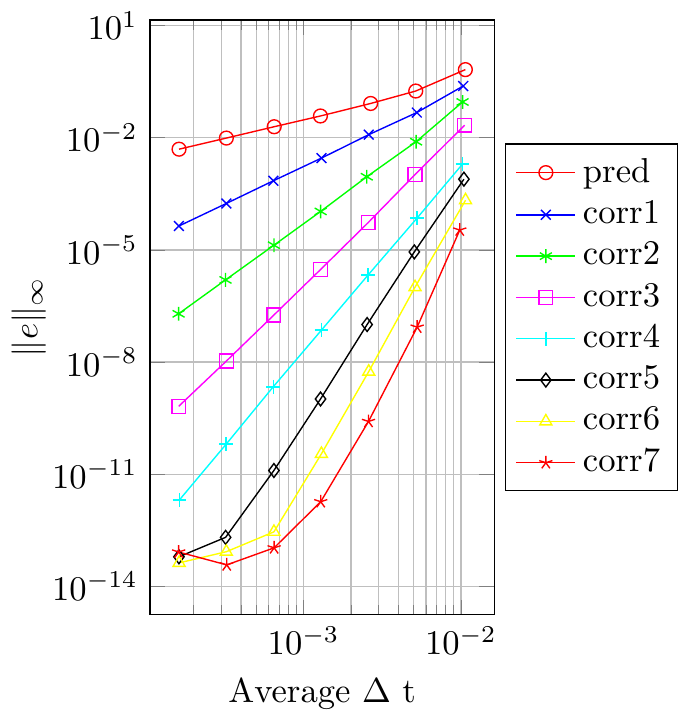}}
  \subfloat[Ratio = 2]{
  \includegraphics{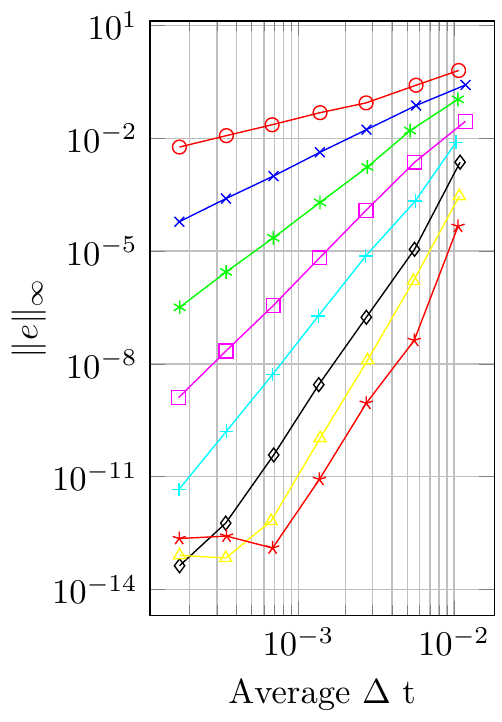}}\\
  \subfloat[Ratio = 4]{
    \includegraphics{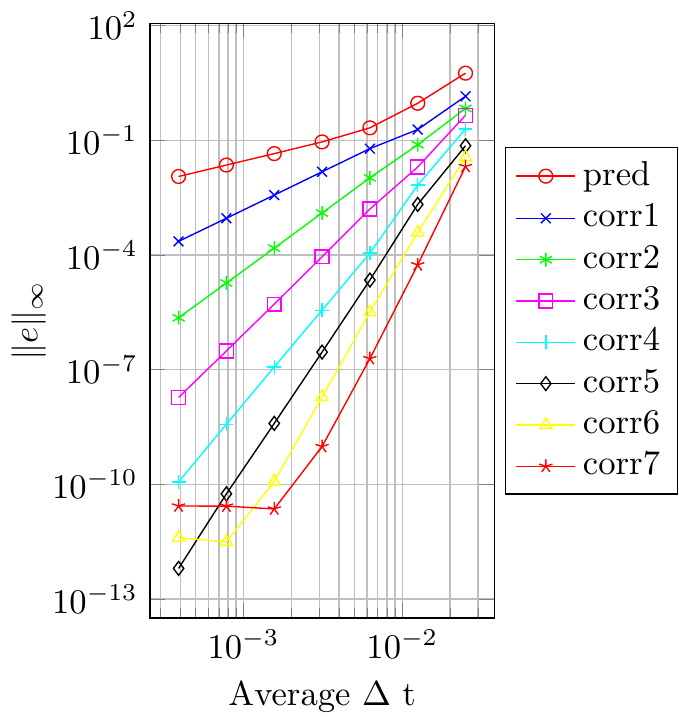}}
  \subfloat[Ratio = 100]{
    \includegraphics{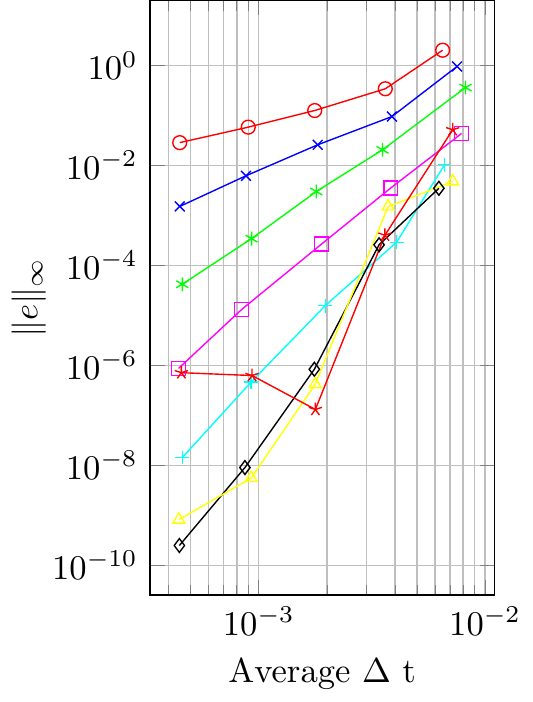}}
  \caption{Lorenz IVP: The design order is illustrated for the RIDC
    methods.  % As before, the performance of higher-order RIDC methods
%     degrade as the ratio of largest to smallest cell increases
  }
  \label{fig:nonunif_lorentz}
\end{figure}

Figure~\ref{fig:nonunif_lorentz} shows the convergence study (error as
a function of mean step-size) for equation~\eqref{eqn:lorenz}. The
reference solution is computed using an RK-45 integrator with a fine
time step. Similar observations can be made that RIDC methods with
non-uniform step-sizes converge with their designed orders of accuracy
(at least up to order 6).

\subsection{Adaptive RIDC}
We study four different variants of RIDC methods with adaptive
step-size control: (i) step-doubling is used for adaptive step-size
control on the prediction level only
(Section~\ref{sec:numerics_predict_double}); (ii) an embedded RK pair is used
for adaptive step-size control on the prediction level only
(Section~\ref{sec:numerics_predict_erk}); (iii) step-doubling is used
for adaptive step-size control on the prediction and correction levels
(Section~\ref{sec:numerics_all_double}); and (iv) step-doubling is
used for adaptive step-size control on the prediction level, and the
computed errors from the error equation \eqref{eqn:error_eqn} are used
for adaptive step-size control on the correction levels.

\subsubsection{Step-Doubling on the Prediction Level Only}
\label{sec:numerics_predict_double}
\begin{figure}[htbp]
  \centering
  \subfloat[prediction]{
    \includegraphics{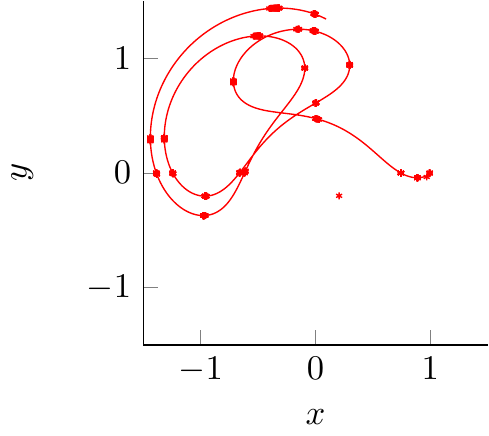}
    \label{fig:orbit_apred_pred}}
  \hfill
  \subfloat[1\textsuperscript{st} Correction]{
    \includegraphics{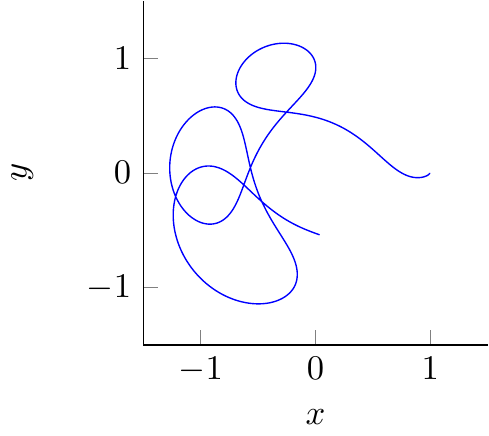}
    \label{fig:orbit_apred_corr1}}\\
  \subfloat[2\textsuperscript{nd} Correction]{
  \includegraphics{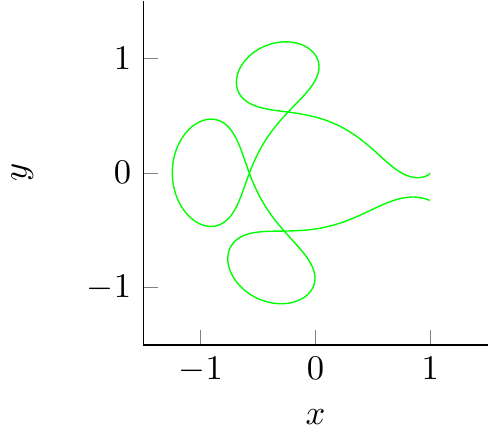}
    \label{fig:orbit_apred_corr2}}
  \hfill
  \subfloat[3\textsuperscript{rd} Correction]{
  \includegraphics{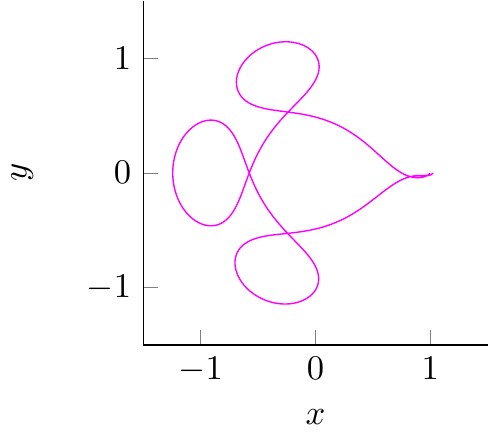}
    \label{fig:orbit_apred_corr3}}
  \caption{Orbit problem: Although the prediction level gives a highly
    inaccurate solution, successive correction loops are able to
    reduce the error and produce the desired orbit.  The red
    circles on the prediction level (a) indicate rejected steps.}
  \label{fig:orbit_apred}
\end{figure}

\begin{figure}[htbp]
  \parbox{0.4\textwidth}{\vspace{0pt}
    \subfloat[Convergence study]{
      \includegraphics{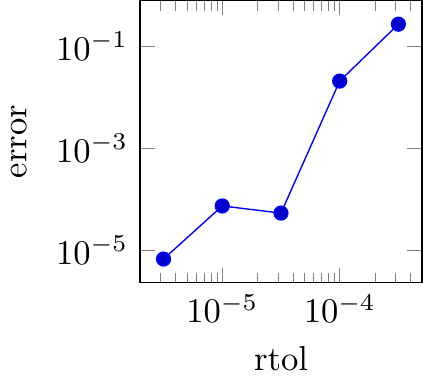}
      \label{fig:orbit_double_conv}
    }
  }
  \begin{tabular}{ccccc}
    \hline
    {\tt rtol} & {\tt atol} & error & naccept & nreject \\\hline
    $10^{-3.5}$ & $10^{-6.5}$ &  2.72e--01    & 1456  & 99  \\
    $10^{-4.0}$ & $10^{-7.0}$ &  2.08e--02    & 2650  & 81  \\
    $10^{-4.5}$ & $10^{-7.5}$ &  5.35e--05    & 4730  & 68  \\
    $10^{-5.0}$ & $10^{-8.0}$ &  7.39e--05    & 8436  & 42  \\
    $10^{-5.5}$ & $10^{-8.5}$ &  6.72e--06    & 15031 & 10  \\
    \hline
  \end{tabular}\\

  \centering
  \subfloat[Adaptive step-sizes selected (reset every 100 steps)]{
    \includegraphics{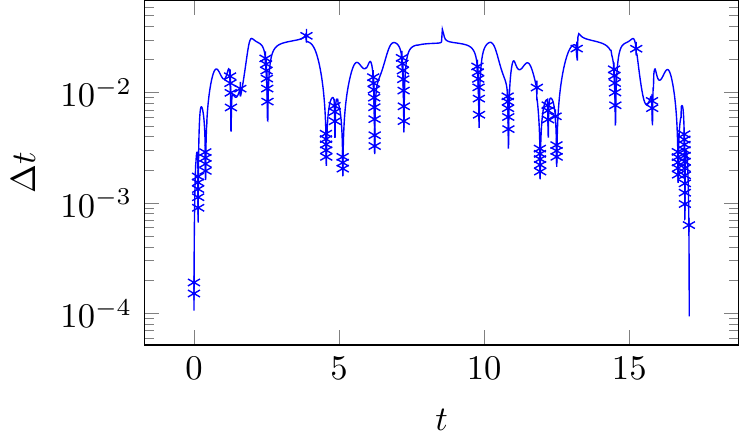}
  \label{fig:orbit_double_100}}
  \quad
  \subfloat[Adaptive step-sizes selected (reset every 400 steps)]{
    \includegraphics{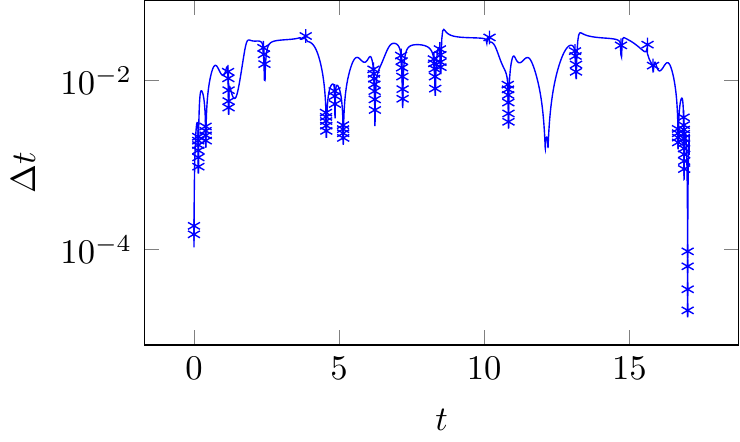}
  \label{fig:orbit_double_400}}

\caption{Orbit Problem: (a) Convergence of a fourth-order RIDC method
  constructed with forward Euler integrators and adaptive step-size
  control on the prediction level (using step-doubling).  Convergence
  is measured relative to the exact solution as the tolerance is
  decreased.  A reset is performed after every 100 accepted steps for
  this convergence study.  In (b), the step-sizes selected for {\tt
    rtol} = $10^{-3.5}$ and {\tt atol} = $10^{-6.5}$ are displayed as
  the solid curve and rejected steps as ``x''s; a reset is performed
  after every 100 steps.  In (c), the reset is performed after every
  400 steps.  Observe that although the number of rejected steps
  increases, the overall $\Delta t$ chosen remains qualitatively
  similar. }
  \label{fig:orbit_apred_conv}
\end{figure}

In this numerical experiment, we solve the orbit
problem~\eqref{eqn:orbit} using a fourth-order RIDC method
(constructed using forward Euler integrators), and adaptive step-size
control on the prediction level only, where step-doubling is used to
provide the error estimate.  As shown in Figure~\ref{fig:orbit_apred},
successive correction loops are able to reduce the error in the
solution and recover the desired orbit.  The red circles in
Figure~\ref{fig:orbit_apred_pred} indicate rejected steps.
Figure~\ref{fig:orbit_double_conv} shows that RIDC with step-doubling
only on the prediction level converges as the tolerance is reduced.
In this experiment, the RIDC integrator is {\em reset} after every 100
accepted steps. By ``reset''~\cite{explicit_ridc}, we mean that the
highest-order solution after every 100 steps is used as an initial
condition to re-initialize the provisional solution; e.g., instead of
solving equation~\eqref{eqn:ode}, one solves a sequence of problems
\begin{align*}
  \left\{\begin{array} {l}
      \displaystyle
      y'(t)  = f(t,y), \quad t \in[t_{100}, \min(b,t_{200})], \\
      \displaystyle
      y(t_{100i}) = \eta_{100(i-1)}^{[P-1]}.
    \end{array}
  \right.
\end{align*}
if $(L-1)$ correctors are applied and $\eta_0^{[L-1]}=y_a$.  The
time steps chosen by the RIDC integrator with resets performed every
100 and 400 steps are shown in Figures~\ref{fig:orbit_double_100} and 
\ref{fig:orbit_double_400}.

In Figure~\ref{fig:orbit_double_100}, $\Delta t_{min} = 1.06 \times
10^{-4}$.  If a non-adaptive fourth-order RIDC method was used with
$\Delta t_{min}$, 160814 uniform time steps would have been required.
By adaptively selecting the time steps for this example and tolerance,
the adaptive RIDC method required approximately one one-hundredth of
the functional evaluations, corresponding to a one hundred-fold
speedup.

\subsubsection{Embedded RK on the Prediction Level Only}
\label{sec:numerics_predict_erk}
\begin{figure}[htbp]
  \parbox{0.4\textwidth}{\vspace{0pt}
    \subfloat[Convergence study]{
      \includegraphics{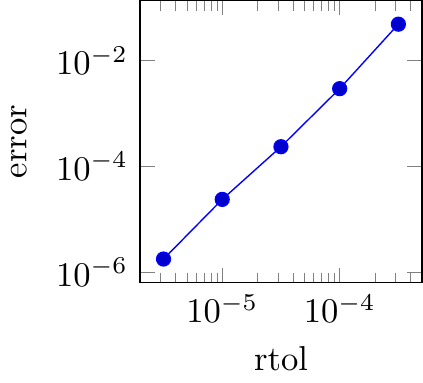}
  \label{fig:orbit_erk_conv}}}
  \begin{tabular}{ccccc}
    \hline
    {\tt rtol} & {\tt atol} & error & naccept & nreject \\\hline
    $10^{-3.5}$ & $10^{-6.5}$ &  4.91e--02    & 2082  & 93  \\
    $10^{-4.0}$ & $10^{-7.0}$ &  2.96e--03    & 3754  & 71  \\
    $10^{-4.5}$ & $10^{-7.5}$ &  2.36e--04    & 6703  & 50  \\
    $10^{-5.0}$ & $10^{-8.0}$ &  2.28e--05    & 11945  & 20  \\
    $10^{-5.5}$ & $10^{-8.5}$ &  1.77e--06    & 21277 & 10  \\
    \hline
  \end{tabular}\\

  \centering
  \subfloat[Adaptive step-sizes selected (reset every 100 steps)]{
    \includegraphics{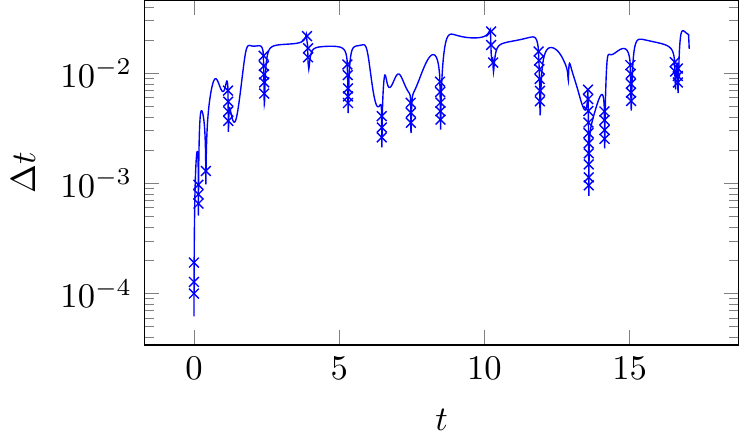}
  \label{fig:orbit_erk_100}}
  \quad
  \subfloat[Adaptive step-sizes selected (reset every 400 steps)]{
    \includegraphics{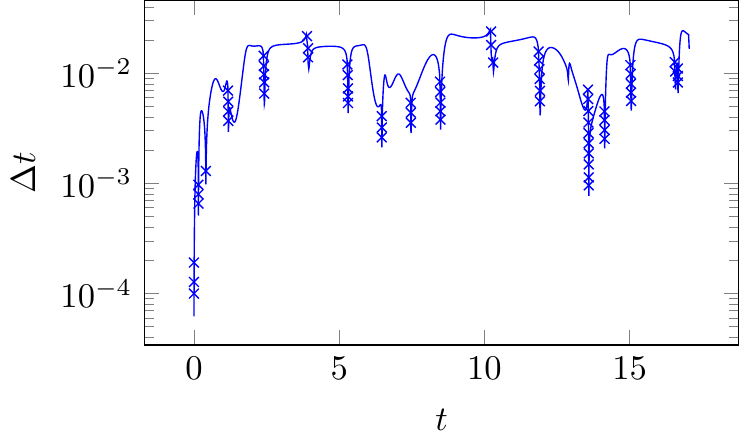}
  \label{fig:orbit_erk_400}}

\caption{Orbit Problem: (a) Convergence of a fourth-order RIDC method
  constructed with forward Euler integrators and adaptive step-size
  control on the prediction level (using an embedded RK pair to
  estimate the error).  Convergence is measured relative to the exact
  solution as the tolerance is decreased.  A reset is performed after
  every 100 accepted steps for this convergence study.  In (b), the
  step-sizes selected for {\tt rtol} = $10^{-3.5}$ and {\tt atol} =
  $10^{-6.5}$ are displayed as the solid curve and rejected steps as
  ``x''s; a reset is performed after every 100 steps.  In (c), the
  reset is performed after every 400 steps.}
  \label{fig:orbit_apred_erk}
\end{figure}

In this numerical experiment, we repeat the orbit
problem~\eqref{eqn:orbit} using a fourth-order RIDC method constructed
again using forward Euler integrators, but the step-size adaptivity on
the prediction level uses a Heun--Euler embedded RK pair. This simple
scheme combines Heun's method, which is second order, with the forward
Euler method,
which is first order.  Figure~\ref{fig:orbit_erk_conv} shows the
convergence of this adaptive RIDC method as the tolerance is reduced.
As the previous example, the RIDC integrator is reset after every 100
accepted steps for the convergence study. In
Figures~\ref{fig:orbit_erk_100} and \ref{fig:orbit_erk_400}, we show
the time steps chosen by the RIDC integrator with resets performed
after 100 or 400 steps, respectively.

\begin{figure}[htbp]
  \begin{minipage}{0.6\textwidth}
    \centering
    \includegraphics{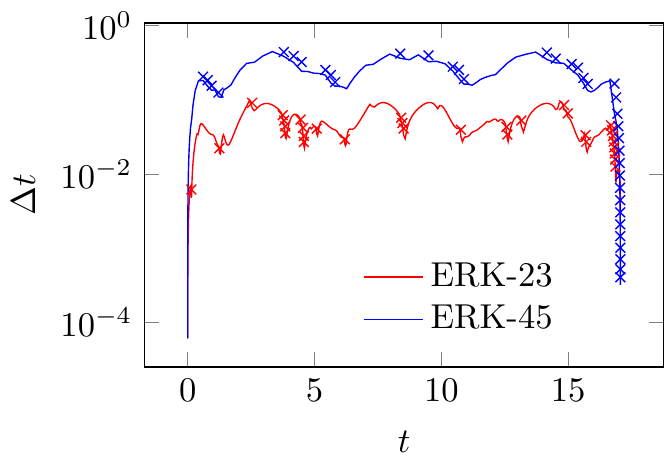}
  \end{minipage}
  \hfill
  \begin{minipage}{0.39\textwidth}
    \caption{\label{fig:orbit_erk} Step-sizes selected by RIDC methods
      constructed using a Bogacki--Shampine method, a 3(2) embedded
      pair (red) and the Runge--Kutta-Fehlberg 4(5) pair.  'x'
      indicates rejected steps.}
  \end{minipage}
\end{figure}

Not surprisingly, the time steps chosen by the RIDC method are
dependent on the specified tolerances and the error estimator (and
consequently the integrators used to obtain a provisional solution to
\eqref{eqn:ode}) used for the control strategy.  One can easily
construct a RIDC integrator using higher-order embedded RK pairs to solve for
a provisional solution to \eqref{eqn:ode}, and then use the forward
Euler method to solve the error equation \eqref{eqn:error_eqn} on subsequent
levels.  For example, Figure~\ref{fig:orbit_erk} shows the step
sizes chosen when the Bogacki--Shampine method \cite{bogacki-shampine}
(a 3(2) embedded RK pair)
% is used to compute the provisional solution (and
%error estimate) for the RIDC integrator.
%Figure~\ref{fig:orbit_erk45_100} shows the step-sizes chosen when the
and the popular Runge--Kutta--Fehlberg 4(5) pair \cite{fehlberg69} is
used to compute the provisional solution (and error estimate) for the
RIDC integrator.  The same tolerance of {\tt rtol} = $10^{-3.5}$ is
used to generate both graphs.  As the order and accuracy of the
predictor increases, one can take larger time steps.  For this
example, using higher-order embedded RK pairs as step-size control
mechanisms for RIDC methods result in less variations in time steps.

\subsubsection{Step-Doubling on All Levels}
\label{sec:numerics_all_double}
As mentioned in Subsection~\ref{sec:adapt_all}, it might be advantageous
to use adaptive step-size control when solving the error equations.
This affords a myriad of parameters that can be used to tune the
step-size control mechanism.  In this set of numerical experiments, we
explore how the choice of tolerances for the prediction/correction
levels affect the step-size selection.

\begin{figure}[htbp]
  \parbox{0.4\textwidth}{\vspace{0pt}
      \includegraphics{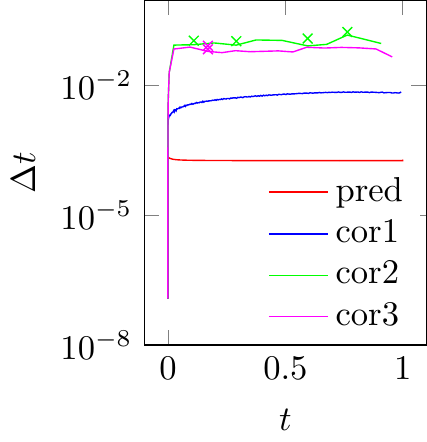}}
  \begin{tabular}{cccccc}
    \hline
    $\ell$ & {\tt rtol} & {\tt atol} & error & naccept & nreject \\\hline
    0 & $10^{-8}$ & $10^{-10}$ & 2.028e--05 & 5480 & 0  \\
    1 & $10^{-8}$ & $10^{-10}$ & 8.824e--07 &  197 & 0  \\
    2 & $10^{-8}$ & $10^{-10}$ & 1.917e--08 &   19 & 4  \\
    3 & $10^{-8}$ & $10^{-10}$ & 6.386e--07 &   25 & 2  \\
    \hline
  \end{tabular}
  \caption{Auzinger IVP: step-size control is implemented on all
    prediction and correction levels.  The same tolerances are used
    for each level.  As expected, the predictor has a hard time
    (forward Euler must satisfy a stringent tolerance); on the other
    hand, life is easy for the correctors.  'x' in the figure
    indicate the rejected steps.}
  \label{fig:auz_aridc_full_dt}
\end{figure}

We first solve the Auzinger IVP using step-doubling on all the levels,
i.e., both predictor and corrector levels.  In
\autoref{fig:auz_aridc_full_dt}, we show the computed step-sizes when
we naively choose the same tolerances on each level.  As expected, the
predictor has to take many steps (to satisfy the stringent
user-supplied tolerance), whereas life is easy for the correctors. In
principle, the correctors are not even needed. Equally important to
note is that the error {\em increases} after the last correction loop.
This might seem surprising at first glance but ultimately may not
unreasonable because the steps selected to solve the third correction
are not based on the solution to the error equation but rather the
original IVP.
% relatively large compared to the steps selected to solve the first
% correction equation.  Large step sizes means that the
% interpolation/quadrature error coefficient is larger, even if the
% interpolant is formally higher order.

\begin{figure}[htbp]
  \parbox{0.4\textwidth}{\vspace{0pt}
    \subfloat[Set 1 of tolerances]{
      \includegraphics{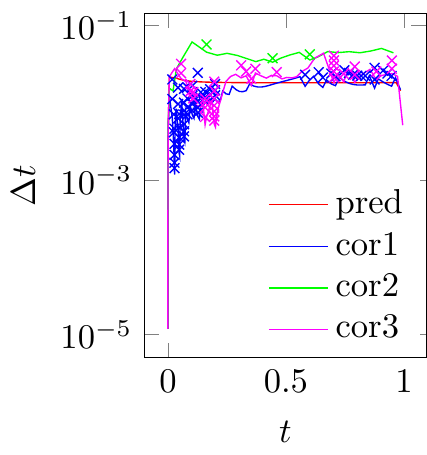}}}
  \begin{tabular}{cccccc}
    \hline
    $\ell$ & {\tt rtol}  & {\tt atol}  & error & naccept & nreject \\\hline
    0 & 1e--04 & 1e--06 & 2.031e--03 & 59 & 0\\
    1 & 1e--06 & 1e--08 & 7.002e--05 & 81 & 61\\
    2 & 1e--08 & 1e--10 & 1.412e--07 & 30 & 3\\
    3 & 1e--10 & 1e--12 & 9.847e--08 & 60 & 33\\\hline
  \end{tabular}
  \parbox{0.4\textwidth}{\vspace{0pt}
    \subfloat[Set 2 of tolerances]{
      \includegraphics{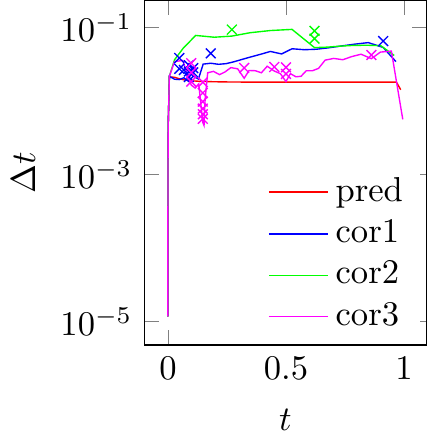}}}
  \begin{tabular}{cccccc}
    \hline
    $\ell$ & {\tt rtol} & {\tt atol} & error & naccept & nreject \\\hline
    0 & 1e--04 & 1e--06 & 2.031e--03 & 59 & 0\\
    1 & 1e--05 & 1e--07 & 1.853e--04 & 31 & 10\\
    2 & 1e--07 & 1e--09 & 1.505e--06 & 21 & 3\\
    3 & 1e--09 & 1e--11 & 9.473e--07 & 44 & 14\\\hline
  \end{tabular}
  \caption{Auzinger IVP: Different tolerances at each level.  With
    the first set of tolerances, the step-size controller for the
    predictor is well behaved, as are the second and third
    correctors.  The step-size controller for the first corrector
    however is noisy.  With the second set of tolerances, the
    step-size controllers for all correctors are reasonably well
    behaved.  }
  \label{fig:auz_aridc_full_dt2}
\end{figure}

Instead of naively choosing the same tolerances on each level, we now
change the tolerance at each level, as described in
\autoref{fig:auz_aridc_full_dt2}.  By making this simple change, the
number of accepted steps on each level are now on the same order of
magnitude.  Not surprisingly, the predictor still selects good
steps. Interestingly in \autoref{fig:auz_aridc_full_dt2}(a), the first
correction is ``noisy'', especially initially.  By picking a different
set of tolerances, we can eliminate the noise, as shown in
\autoref{fig:auz_aridc_full_dt2}(b).

\subsubsection{Using Solutions from the Error Equation}
\label{sec:numerics_successive}

\begin{figure}[htbp]
  \parbox{0.4\textwidth}{\vspace{0pt}
    \includegraphics{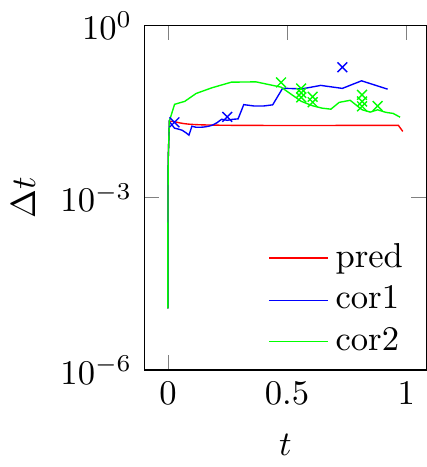}}
  \begin{tabular}{cccccc}
    \hline
    $\ell$ & {\tt rtol} & {\tt atol} & error & naccept & nreject \\\hline
    0 & 1e--04 & 1e--06 & 2.031e--03 & 59 & 0\\
    1 & 1e--04 & 1e--06 & 7.249e--04 & 33 & 3 \\
    2 & 1e--04 & 1e--06 & 6.513e--06 & 26 & 10
  \end{tabular}
  \parbox{0.4\textwidth}{\vspace{0pt}
    \includegraphics{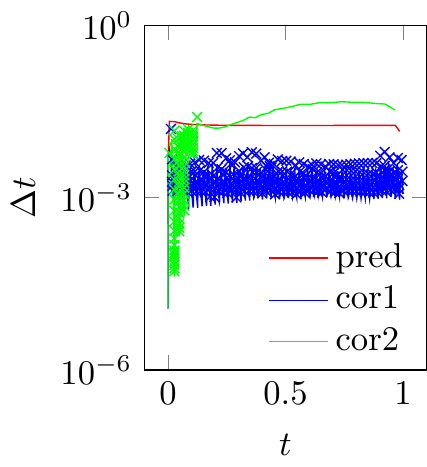}}
  \begin{tabular}{cccccc}
    \hline
    $\ell$ & {\tt rtol} & {\tt atol} & error & naccept & nreject \\\hline
    0 & 1e--04 & 1e--06 & 2.031e--03 & 59 & 0\\
    1 & 1e--05 & 1e--07 & 1.063e--05 & 657 & 305 \\
    2 & 1e--06 & 1e--08 & 9.446e--08 & 75 & 76
  \end{tabular}  \parbox{0.4\textwidth}{\vspace{0pt}
    \includegraphics{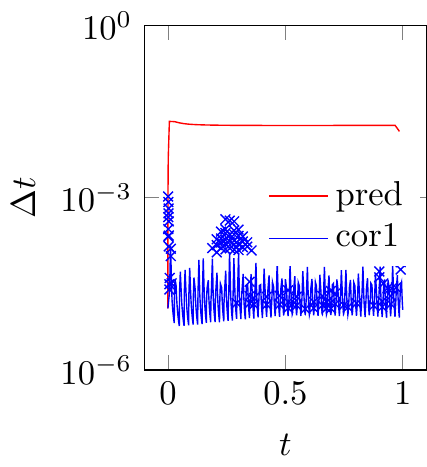}}
  \begin{tabular}{cccccc}
    \hline
    $\ell$ & {\tt rtol} & {\tt atol} & error & naccept & nreject \\\hline
    0 & 1e--04 & 1e--06 & 2.031e--03 & 59 & 0\\
    1 & 1e--07 & 1e--09 & 1.178e--07 & 60571 & 94
  \end{tabular}
  %% \parbox{0.4\textwidth}{\vspace{0pt}
  %%   \includegraphics[width=0.4\textwidth]{figures/auz_aridc_dc_dt4}}
  %% \begin{tabular}{cccccc}
  %%   \hline
  %%   $\ell$ & {\tt rtol} & {\tt atol} & error & naccept & nreject \\\hline
  %%   0 & 1e--06 & 1e--08 & 2.028e--04 & 552 & 0\\
  %%   1 & 1e--07 & 1e--09 & 1.074e--07 & 6561 & 3081
  %% \end{tabular}
  \caption{Auzinger Problem: Step-doubling on prediction levelz, using
    successive levels for error estimation for step control on the
    error equation.  Step-size controller for the corrector is 
    noisy.}  
  \label{fig:auzinger_formulation3}
\end{figure}

\begin{figure}
  \parbox{0.4\textwidth}{\vspace{0pt}
    \includegraphics{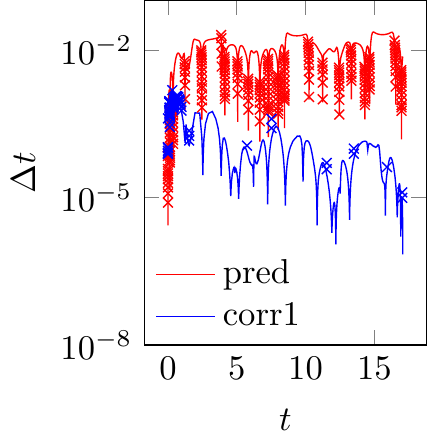}}
  \begin{tabular}{cccccc}
    \hline
    $\ell$ & {\tt rtol} & {\tt atol} & error & naccept & nreject \\\hline
    0 & 1e--04 & 1e--04 & 2.405e--00 & 2261 & 230\\
    1 & 1e--04 & 1e--04 & 7.234e--01 & 475181 & 84
  \end{tabular}
  \caption{Orbit Problem: Step-doubling on prediction level, using
    successive levels for error estimation for step control on the
    error equation. }
  \label{fig:orbit_formulation3}
\end{figure}

As mentioned in Subsection~\ref{sec:using_error_equations}, using the
solution from the error equation~\eqref{eqn:error_eqn} as the local error
estimate for step-size control on a given level is potentially
problematic because the step-size controller can only control the local
error introduced on that level whereas the true local error generally
contains contributions from all previous levels. For completeness, we
present the results of this adaptive RIDC formulation applied to the
Auzinger problem (\autoref{fig:auzinger_formulation3}) and the Orbit
problem (\autoref{fig:orbit_formulation3}).  Step-doubling is used for
step-size adaptivity on the predictor level, solutions from the error
equation are used to control step-sizes for the corrector levels.  For
the Auzinger problem, we observe in the top figure that if the
tolerances are held fixed on each level, each correction level
improves the solution.  If the tolerance is reduced slightly on each
level, the step-size controller gives a poor step-size selection (many
rejected steps), even for this smoothly varying problem.  For the
Orbit IVP, \autoref{fig:orbit_formulation3} shows that the corrector
improves the solution if the tolerances are held fixed at all levels;
however the corrector requires {\em many} steps.  A second correction
loop was not attempted. Reducing the tolerance for the first corrector
resulted in inordinately many rejected steps.

\section{Conclusions}
\label{sec:conclusions}
In this paper, we formulated RIDC methods that incorporate local error
estimation and adaptive step-size control.  Several formulations were
discussed in detail: (i) step-doubling on the prediction level, (ii)
embedded RK pairs on the prediction level, (iii) step-doubling on the
prediction and error levels, and (iv) step-doubling for the prediction
level but using the solution from the error equation for step-size
control; other formulations are also alluded to.  A convergence
theorem from~\cite{XXS} can be extended to RIDC methods that use
adaptive step-size control on the prediction level.  Numerical
experiments demonstrate that RIDC methods with non-uniform steps
converge as designed and illustrate the type of behavior that might be
observed when adaptive step-size control is used on the prediction and
correction levels.  Based on our numerical study, we conclude that
adaptive step-size control on the prediction level is viable for RIDC
methods.  In a practical application where a user gives a specified
tolerance, this prescribed tolerance must be transformed to a specific
tolerance that is fed to the predictor.

%1) another area is for prediction only error control: user gives
%tolerance.  We have to choose an practical tolerance to feed to the
%prediction (e.g., user gives $10^{-12}$ we use $10^{-12/p}$

%2) same, but for multi-level error control: generate tolerances for
%each level.

\section*{Acknowledgments}
This publication was based on work supported in part by Award No
KUK-C1-013-04, made by King Abdullah University of Science and
Technology (KAUST), AFRL and AFOSR under contract and grants
FA9550-12-1-0455, NSF grant number DMS-0934568, NSERC grant number
RGPIN-228090-2013, and the Oxford Center for Collaborative and Applied
Mathematics (OCCAM).

\bibliography{ridc}
\bibliographystyle{ieeetr}

\end{document}